\documentclass[a4paper, 12pt]{article}
\usepackage{mathtext} 
\usepackage[T2A, T1]{fontenc}
\usepackage[utf8x]{inputenc}
\usepackage[english, polish, russian]{babel}
\usepackage[fleqn,tbtags]{amsmath}
\usepackage{mathtools}
\usepackage[left=1cm, right=1cm, top=1cm, bottom=2cm, bindingoffset=0cm]{geometry}
\usepackage{amsfonts}
\usepackage{amssymb,amsmath}
\usepackage{amsthm}
\usepackage{tikz}
\usepackage{tikz-cd}
\usepackage{stackrel}
\usepackage{soul}
\usepackage{enumitem} 
\usepackage[normalem]{ulem}
\usepackage{etoolbox}
\patchcmd{\thebibliography}{\section*{\refname}}{}{}{}

\newtheorem{theorem}{Теорема}
\newtheorem*{lemma*}{Лемма}

\newenvironment{definition}[1][Определение.]{\begin{trivlist}
\item[\hskip \labelsep {\bfseries #1}]}{\end{trivlist}}

\newenvironment{remark}[1][Замечание.]{\begin{trivlist}
\item[\hskip \labelsep {\bfseries #1}]}{\end{trivlist}}

\DeclarePairedDelimiterX\Set[2]{\lbrace}{\rbrace}%
 { #1 \,\delimsize|\, #2 }

\DeclarePairedDelimiter\abs{\lvert}{\rvert}

\DeclarePairedDelimiter\norm{\lVert}{\rVert}

\newcommand\restr[2]{{
  \left.\kern-\nulldelimiterspace 
  #1
  \vphantom{\big|} 
  \right|_{#2} 
}}

\newcommand{\genlb}[1]{%
\begin{tikzpicture}[#1]%
\draw (0.5,0) -- (0.5,1);%
\draw (0,0.5) -- (1,0.5);%
\path [fill=black] (0, 0) rectangle (0.5, 0.5);%
\end{tikzpicture}%
}
\newcommand{\genlt}[1]{%
\begin{tikzpicture}[#1]%
\draw (0.5,0) -- (0.5,1);%
\draw (0,0.5) -- (1,0.5);%
\path [fill=black] (0, 0.5) rectangle (0.5, 1);%
\end{tikzpicture}%
}
\newcommand{\genrt}[1]{%
\begin{tikzpicture}[#1]%
\draw (0.5,0) -- (0.5,1);%
\draw (0,0.5) -- (1,0.5);%
\path [fill=black] (0.5, 0.5) rectangle (1, 1);%
\end{tikzpicture}%
}
\newcommand{\genrb}[1]{%
\begin{tikzpicture}[#1]%
\draw (0.5,0) -- (0.5,1);%
\draw (0,0.5) -- (1,0.5);%
\path [fill=black] (0.5, 0) rectangle (1, 0.5);%
\end{tikzpicture}%
}

\newcommand{\genl}[1]{%
\begin{tikzpicture}[#1]%
\draw (0.5,0) -- (0.5,1);%
\draw (0,0.5) -- (1,0.5);%
\path [fill=black] (0, 0) rectangle (0.5, 1);%
\end{tikzpicture}%
}
\newcommand{\gent}[1]{%
\begin{tikzpicture}[#1]%
\draw (0.5,0) -- (0.5,1);%
\draw (0,0.5) -- (1,0.5);%
\path [fill=black] (0, 0.5) rectangle (1, 1);%
\end{tikzpicture}%
}
\newcommand{\genr}[1]{%
\begin{tikzpicture}[#1]%
\draw (0.5,0) -- (0.5,1);%
\draw (0,0.5) -- (1,0.5);%
\path [fill=black] (0.5, 0) rectangle (1, 1);%
\end{tikzpicture}%
}
\newcommand{\genb}[1]{%
\begin{tikzpicture}[#1]%
\draw (0.5,0) -- (0.5,1);%
\draw (0,0.5) -- (1,0.5);%
\path [fill=black] (0, 0) rectangle (1, 0.5);%
\end{tikzpicture}%
}

\newcommand{\gennlb}[1]{%
\begin{tikzpicture}[#1]%
\draw (0.5,0) -- (0.5,1);%
\draw (0,0.5) -- (1,0.5);%

\path [fill=black] (0, 0.5) rectangle (1, 1);%
\path [fill=black] (0.5, 0) rectangle (1, 1);%
\end{tikzpicture}%
}
\newcommand{\gennlt}[1]{%
\begin{tikzpicture}[#1]%
\draw (0.5,0) -- (0.5,1);%
\draw (0,0.5) -- (1,0.5);%

\path [fill=black] (0.5, 0) rectangle (1, 1);%
\path [fill=black] (0, 0) rectangle (1, 0.5);%
\end{tikzpicture}%
}
\newcommand{\gennrt}[1]{%
\begin{tikzpicture}[#1]%
\draw (0.5,0) -- (0.5,1);%
\draw (0,0.5) -- (1,0.5);%

\path [fill=black] (0, 0) rectangle (0.5, 1);%
\path [fill=black] (0, 0) rectangle (1, 0.5);%
\end{tikzpicture}%
}
\newcommand{\gennrb}[1]{%
\begin{tikzpicture}[#1]%
\draw (0.5,0) -- (0.5,1);%
\draw (0,0.5) -- (1,0.5);%

\path [fill=black] (0, 0) rectangle (0.5, 1);%
\path [fill=black] (0, 0.5) rectangle (1, 1);%
\end{tikzpicture}%
}

\newcommand{\rbs}{\genrb{scale=0.25}}

\newcommand{\ts}{\gent{scale=0.25}}
\newcommand{\rs}{\genr{scale=0.25}}
\newcommand{\bs}{\genb{scale=0.25}}

\newcommand{\nlbs}{\gennlb{scale=0.25}}

\title{K-замкнутость весовых пространств Харди на $\mathbb{T}^2$}
\date{2017}
\author{Боровицкий Вячеслав Андреевич}

\begin{document}

\maketitle

\vspace*{4cm}

\begin{abstract}
В работе доказывается, при некоторых условиях на веса, K-замкнутость пары весовых пространств Харди на двумерном торе в паре соответствующих весовых пространств Лебега. Вопрос K-замкнутости пространств Харди на двумерном торе до этого рассматривался лишь в безвесовом случае, либо для весов, распадающихся в произведение двух функций одной переменной (так называемых "разделяющихся весов"). В работе рассмотрен случай некоторых неразделяющихся весов.
\end{abstract}

\newpage

					\section*{\begin{center}\S1 Введение\end{center}}

					\hspace{\parindent} В данной работе мы рассматриваем некоторые вопросы о $K$-замкнутости пары весовых пространств Харди на двумерном торе $H^r(w_1(\cdot, \cdot)), H^p(w_2(\cdot, \cdot))$ в паре соответствующих им весовых пространств Лебега $L^r(w_1(\cdot, \cdot)), L^p(w_2(\cdot, \cdot))$. \\

					Пусть $(X_1, X_2)$ совместимая пара банаховых или квази-банаховых пространств (то есть они вложены в некоторое объемлющее топологическое векторное пространство), $Y_1$ и $Y_2$ --- замкнутые подпространства соответственно в $X_1$ и $X_2$. Введем определение.

					\begin{definition}
					Пара $(Y_1, Y_2)$ называется $K$-замкнутой в паре $(X_1, X_2)$, если существуют такие две абсолютные константы $C_1,C_2$, что для всех элементов $f \in Y_1 + Y_2$, $g \in X_1$, $h \in X_2$ таких, что $f=g+h$, найдутся такие элементы $g' \in Y_1$, $h' \in Y_2$, что $f=g'+h'$ и при этом $\norm*{g'}_{X_1} \leq C_1 \norm*{g}_{X_1}$, $\norm*{h'}_{X_2} \leq C_2 \norm*{h}_{X_2}$.
					\end{definition}

					Теоремы о $K$-замкнутости интересны сами по себе, так как в таком случае пара $(Y_1, Y_2)$ наследует многие интерполяционные свойства пары $(X_1, X_2)$, в частности верна формула $(Y_1, Y_2)_{\theta, q} = (X_1, X_2)_{\theta, q} \cap (Y_1 + Y_2)$ для интерполяции пространств вещественным методом. \\

					Напомним, что классические пространства Харди на $n$-мерном торе это
					$$H^p(\mathbb{T}^n) = \text{Clos} \ \text{Lin} \Set*{z_1^{j_1} z_2^{j_2} \dots z_n^{j_n}}{j_1, j_2, \dots, j_n \in \mathbb{N} \cup \{0\}},$$
					где $0 < p \leq \infty$, при $p < \infty$ замыкание берется в $L^p(\mathbb{T}^n)$-норме, а при $p = \infty$ в *-слабой топологии. \\

					Вопрос $K$-замкнутости пары $(H^r(\mathbb{T}^n), H^p(\mathbb{T}^n))$ в паре $(L^r(\mathbb{T}^n), L^p(\mathbb{T}^n))$ для $n=1$ решен давно и может уже считаться классическим, причем исторически интерполяционные свойства шкалы $H^p(\mathbb{T})$ изучались еще до появления понятия $K$-замкнутости (по поводу истории вопроса см. \cite{kis2}). Рассмотрение безвесового случая при $n=2$ было завершено несколько позже, статьей \cite{KisXu}. В ней был решен самый сложный случай $p = \infty$. В случае $n \geq 3$ ``мягкими'' методами решен вопрос $K$-замкнутости для $r,p < \infty$ (см. \cite{xu}), а вот при $p = \infty$ ничего не понятно уже даже для $n=3$.\\

					Перед тем как начать говорить о проблемах $K$-замкнутости весовых пространств Харди, напомним еще одно определение. Мы говорим, что вес $w: \mathbb{T} \to (0, +\infty)$ удовлетворяет условию Макенхаупта $A_p$, $1 < p < \infty$, тогда и только тогда, когда  
					$$\left(\frac{1}{|B|}\int_B w(t)dt\right) \left(\frac{1}{|B|}\int_B w(t)^{\frac{1}{1-p}}dt\right)^{p-1} \leq C < \infty$$
					для любой дуги $B$. Вес $w$ удовлетворяет условию $A_1$ тогда и только тогда, когда
					$$\frac{1}{|B|}\int_B w(t)dt \leq C w(x),$$
					для любой дуги $B$ и для любого $x \in B$. Вес $w$ удовлетворяет $A_\infty$ тогда и только тогда, когда $w$ удовлетворяет какому-нибудь $A_p$ для $1 \leq p < \infty$. При этом наилучшие константы $C$, фигурирующие в правых частях неравенств, принято называть $A_p$-константами, а вместо ``$w$ удовлетворяет условию $A_p$'' иногда пишут $w \in A_p$. Когда мы говорим о весах двух переменных $w: \mathbb{T}^2 \to (0, \infty)$, условия $A_p$ определяются таким же образом, только вместо дуг символ $B$ будет обозначать множества вида $B_1 \times B_2$, где $B_1, B_2$ --- дуги. \\

					Некоторые весовые результаты о $K$-замкнутости пространств Харди в пространствах Лебега могут быть легко получены из соответствующих безвесовых результатов, наложением на вес соответствующих ситуации условий Макенхаупта $A_p$, но такие условия часто являются слишком ограничительными. \\

					Толчком к изучению интерполяции (в основном в форме $K$-замкнутости) весовых пространств Харди при минимальных ограничениях на веса в случае одной переменной стало возникновение этого вопроса в контексте доказательства теоремы Гротендика для пространства, сопряженного к диск-алгебре (диск-алгебра --- это замыкание аналитических полиномов в норме пространства непрерывных функций на окружности), см. подробнее в обзоре \cite{KisAbsSumm}. В результате такие условия были найдены. Верна следующая теорема.
						\begin{theorem}
							\label{super_old}
							Пусть $w_1, w_2: \mathbb{T} \to (0, +\infty)$ --- веса, для которых $\log w_1(\cdot) \in L^1(\mathbb{T})$,  $\log w_2(\cdot) \in L^1(\mathbb{T})$. Пусть $0 < r < p \leq \infty$. Тогда пара $(H^r(w_1(\cdot)),H^p(w_2(\cdot)))$ $K$-замкнута в паре $(L^r(w_1(\cdot)),L^p(w_2(\cdot)))$ тогда и только тогда, когда $\log \frac{w_1^{1/r}(\cdot)}{w_2^{1/p}(\cdot)} \in BMO(\mathbb{T})$ (при $p=\infty$ условие $\log w_1^{1/r}(\cdot) w_2(\cdot) \in BMO(\mathbb{T})$).
						\end{theorem}

						Отметим, что условия $\log w_1(\cdot) \in L^1(\mathbb{T})$, $\log w_2(\cdot) \in L^1(\mathbb{T})$ являются скорее не ограничением применимости теоремы, а очерчивают обстоятельства, при которых не происходит вырождения в определении весовых пространств Харди. Кроме того, в случае $r=p$ теорема остается справедливой. \\
						
						Остановимся ненадолго на этом вопросе. До настоящего момента мы так и не ввели определение весовых пространств Харди. Вообще говоря, это довольно сложный и тонкий вопрос: например для того, чтобы дать определение, аналогичное данному выше для безвесовых пространств Харди, нам нужно было бы, чтобы аналитические полиномы, чье замыкание мы берем, лежали в $H^p(w)$, для чего нужно было бы требовать $w \in L^1(\mathbb{T}^n)$, а это может быть слишком ограничительным условием. Но в случае одной переменной этот вопрос решается довольно простым и естественным образом. Для веса $w: \mathbb{T} \to (0, \infty)$ с суммируемым логарифмом (а мы, как уже говорилось выше, рассматриваем только такие веса), находим внешнюю функцию $u$ такую, что $\abs*{u}=w$; тогда можем определить $H^p(w) = \Set*{f / u^{\frac{1}{p}}}{f \in H^p(\mathbb{T})}$ с нормой $\norm*{g}_{H^p(w)} = \norm*{g u^{\frac{1}{p}}}_{H^p(\mathbb{T})}$. В случае $n = 2$, который будет рассматриваться в дальнейшем, при $0<p<\infty$ мы будем использовать следующее определение: будем рассматривать веса $w: \mathbb{T}^2 \to (0, \infty)$ вида $w(z_1, z_2) = a(z_1) u(z_1, z_2) b(z_2)$, где $u \in L^1(\mathbb{T}^2)$, а у функций $a$ и $b$ суммируемые логарифмы (пусть $\widetilde{a}$, $\widetilde{b}$ внешние функции, построенные по $a$ и $b$), тогда $H^p(w) = \Set*{f / (\widetilde{a} \widetilde{b})^{\frac{1}{p}}}{f \in H^p(u(\cdot, \cdot))}$ с нормой $\norm*{g}_{H^p(w)} = \norm*{g (a b)^{\frac{1}{p}}}_{H^p(u(\cdot, \cdot))}$, где $H^p(u(\cdot, \cdot))$ определяется как замыкание аналитических полиномов в $L^p(u(\cdot, \cdot))$-норме, взятое с этой же нормой. Для произвольного веса (с суммируемым логарифмом) $w: \mathbb{T}^2 \to (0, \infty)$, пространство $L^\infty(w(\cdot, \cdot))$ мы понимаем как
					$$L^\infty(w(\cdot, \cdot)) := \Set*{f:\mathbb{T}^2 \to \mathbb{C}}{\text{ess sup} \Set*{f(z_1, z_2)/w(z_1, z_2)}{(z_1, z_2) \in \mathbb{T}^2} < \infty}$$
					с естественной нормой. Пространство $H^\infty(w(\cdot, \cdot))$ для веса $w$ вида $w(z_1, z_2) = a(z_1) u(z_1, z_2) b(z_2)$, где $u \in L^1(\mathbb{T}^2)$, а у функций $a$ и $b$ суммируемые логарифмы, определяем как аннулятор пространства $L^P_1(w)$, которое, в свою очередь, определяется аналогично $H_1(w)$, только с заменой аналитических многочленов на многочлены со спектром во множестве $\mathbb{N} \times \mathbb{Z} \cup \mathbb{Z} \times \mathbb{N}$, (множество функций со спектром в таком множестве мы будем обозначать $\nlbs$). Отметим, что трюк с внешними функциями в данной ситуации также применим, так как принадлежность функции $\nlbs$ инвариантна относительно домножения на аналитические функции. Двойственность, которую мы только что использовали и будем продолжать использовать всегда в дальнейшем, это $\int f \bar{g}$. При таком определении двойственности и пространства $L^\infty(w(\cdot, \cdot))$ верно $L^1(w(\cdot, \cdot))^* = L^\infty(w(\cdot, \cdot))$. \\

					Для интерполяции пространств Харди двух переменных к настоящему времени была установлена следующая теорема (она доказана в \cite{kis}).\\

						\begin{theorem}
							\label{main_old}
							Пусть $1 < r < \infty$; $w_1(z_1, z_2) = a_1(z_1) b_1(z_2)$, $w_2(z_1, z_2) = a_2(z_1) b_2(z_2)$, где функции $a_i$, $b_i$ удовлетворяют условию $\log a_1(\cdot), \log b_1(\cdot), \log a_2(\cdot), \log b_2(\cdot) \in BMO(\mathbb{T})$. \\
							Тогда пара $(H^r(w_1(\cdot, \cdot)),H^\infty(w_2(\cdot, \cdot)))$ $K$-замкнута в паре $(L^r(w_1(\cdot, \cdot)),L^\infty(w_2(\cdot, \cdot)))$.
						\end{theorem}
						Отметим, что, на самом деле, в оригинальной теореме вместо условия на принадлежность самих логарифмов весов к $BMO(\mathbb{T})$, фигурируют условия с логарифмами отношений весов, аналогичные условиям из теоремы \ref{super_old}. Таким образом аналог одномерного условия оказывается достаточным, когда вес $w(\cdot, \cdot)$ разделяется в произведение двух весов от одной переменной. \\

						В данной работе мы, во-первых, хотели показать, что при некоторых условиях на функцию $u(\cdot, \cdot)$, теорема \ref{main_old} верна для весов вида $w(z_1, z_2) = a(z_1)u(z_1, z_2)b(z_1)$. Наш основной результат на эту тему изложен в следующей теореме.

\begin{theorem}
\label{rght}
$K$-замкнутость пары
$$\left(H_p(a_1(z_1) u_1(z_1, z_2) b_1(z_2)), H_\infty(a_2(z_1) u_2(z_1,z_2) b_2(z_2))\right)$$
в паре
$$\left(L_p(a_1(z_1) u_1(z_1, z_2) b_1(z_2)), L_\infty(a_2(z_1) u_2(z_1,z_2) b_2(z_2))\right)$$
имеет место при следующих условиях:
	\begin{enumerate}[label={\arabic*)}]
		\item $u_1$ удовлетворяет двумерному $A_p$,
		\item $u_2$ удовлетворяет двумерному $A_1$,
		\item $\log(a_i), \log(b_i) \in \text{BMO}$,
		\item $u_2^p u_1$ удовлетворяет $A_{\infty}$ по второй переменной равномерно (в смысле равномерной константы в обратном неравенстве Гельдера),
	\end{enumerate}
\end{theorem}

Заметим, что вышеописанная теорема несколько отличается от представленной в соответствующем разделе. Мы специально разместили здесь версию с более простыми условиями, хотя и не в максимальной доказанной нами общности. Вывести теорему \ref{rght} из теоремы \ref{inf_neib_all_q} читатель сможет без труда. \\

						Во-вторых, соединив идеи из \cite{xu} и \cite{uncond}, мы доказали следующее весовое утверждение для случая $n=2$ и $0 < r \leq 1 < p < \infty$.

\begin{theorem}
	\label{lft}
	Если веса $w_1(\cdot, \cdot)$, $w_2(\cdot, \cdot)$ удовлетворяют условиям
	$$w_1(\cdot, \cdot) \in A_\infty,\qquad w_2(\cdot, \cdot) \in A_p,$$
	то пара $\left(H^r(w_1(\cdot, \cdot)), H^p(w_2(\cdot, \cdot))\right)$ K-замкнута в $\left(L^r(w_1(\cdot, \cdot)), L^p(w_2(\cdot, \cdot))\right)$.
\end{theorem}
						Заметим опять, что в соответствующем разделе доказана несколько более сильная теорема, а вышеприведенный результат выбран как основной, поскольку его легче использовать. Кроме того, похожую теорему, по-видимому, можно сформулировать и для любой размерности $n$, но в данной работе мы ограничились случаем $n=2$. \\

						Наконец, используя вместе предыдущие два результата, мы получили методом ``склейки трех шкал'' следующую теорему, аналогичную общему одномерному результату о $K$-замкнутости подпространств, порожденных сингулярными интегральными проекторами, сформулированному в \cite{ruts}.
\begin{theorem}
	\label{glue}
	Если $w_1, w_2 \in A_1$ и $w_1 w_2 \in A_\infty$, то пара 
	$$\left(H_1(w_1(z_1, z_2)), H_\infty(w_2(z_1,z_2))\right)$$
	$K$-замкнута в паре
	$$\left(L_1(w_1(z_1, z_2)), L_\infty(w_2(z_1,z_2))\right).$$
\end{theorem}

						В работе, помимо введения, три параграфа. \S2 содержит результаты, касающиеся $K$-замкнутости в окрестности единицы: то есть при $0 < r \leq 1 < p < \infty$. В нем будет доказана теорема \ref{lft}. В \S3 рассматривается случай $1<r<p=\infty$, доказывается теорема \ref{rght}. В \S4 из результатов \S2 и \S3 выводится теорема \ref{glue}, про случай $r=1, p=\infty$. \\

						Добавим еще пару слов об обозначениях, которые проходят сквозь всю работу. \\

						В дальнейшем нам несколько раз нужно будет говорить о классах функций со спектром в определенном координатном угле, поэтому мы вводим для них специальные символы, а именно: $$\genlb{scale=0.3}, \genlt{scale=0.3}, \genrt{scale=0.3}, \genrb{scale=0.3}, \gent{scale=0.3}, \genr{scale=0.3}, \genl{scale=0.3}, \genb{scale=0.3}, \gennlb{scale=0.3}, \gennlt{scale=0.3}, \gennrt{scale=0.3}, \gennrb{scale=0.3},$$ которые читаются так: функции, спектр которых лежит в закрашенной черным области (черные линии, где видны, символизируют собой координатные оси, они добавлены для удобства). Ось абсцисс будет всегда отвечать первой переменной, ось ординат --- второй. Включать или не включать функции со спектром на границе будет ясно из контекста, в котором символ используется. \\

						Для того, чтобы несколько уменьшить громоздкость выкладок, мы вводим символы $\lesssim$ и $\gtrsim$, которые будут обозначать оценки с константами, значение которых нас не интересует. То есть обозначение $A \lesssim B$ значит, что выполнено $A \leq C B$ для некоторой положительной константы $С$. \\

						Символом $\mu$ будем обозначать нормированную меру Лебега на одномерном торе (то есть такую, что $\mu(\mathbb{T}) = 1$). Для $X \subseteq \mathbb{T}$ запись $(w \mu)(X)$ будет обозначать $\int_X w$, то есть весовую меру множества $X$. \\
						
						Кроме этого, нам будут встречаться (на самом деле однажды уже встречались) обозначения вида $X^Q$, где $X$ --- какая-то квази-банахова решетка измеримых функций, а $Q$ некоторый проектор. При этом $Q$ не обязан действовать в пространстве $X$: если он все же действует в $X$, то, стандартным образом, $X^Q = \Set*{f \in X}{Qf=f}$, если же нет, то мы будем фиксировать за проектором $Q$ некоторое линейное подпространство $D \subseteq X$ (в реальных случаях оно, чаще всего, будет плотным), на котором $Q$ определен и принимает значения в $X$, а пространство $X^Q$ определять как $Clos \Set*{f \in D}{Qf=f}$ (то, что $Q$ проектор, в последнем случае понимается как $Q^2=Q$ на множестве $D$). Введем сразу важнейший для \S3 оператор $P$. Это проектор на $\nlbs$; более точно, он действует на тригонометрические многочлены двух переменных, зануляя коэффициенты при степенях в множестве $(\mathbb{Z}\setminus \mathbb{N}) \times (\mathbb{Z}\setminus \mathbb{N})$.
						Теперь, следуя вышепреведенной общей конструкции, можем определить решетку $L_s^P(u(\cdot, \cdot))$, где $u \in L^1(\mathbb{T}^2)$, $s < \infty$ (в роли множества $D$, как и всегда для проектора $P$, будет выступать множество тригонометрических многочленов). Но все же таким образом мы не можем определить $L_s^P(w(\cdot, \cdot))$, где $w$ вида $w(z_1, z_2) = a(z_1) u(z_1, z_2) b(z_2)$, $s < \infty$, функция $u \in L^1(\mathbb{T}^2)$, а у функций $a$ и $b$ суммируемые логарифмы. Случай $s=\infty$, исключенный здесь нами, получается из соображений двойственности, как это уже было сделано на несколько абзацев выше, в этом абзаце мы более не будем его обсуждать. Чтобы справиться с весами вида $w(z_1, z_2) = a(z_1) u(z_1, z_2) b(z_2)$, мы определяем пространство $L_s^P(w(\cdot, \cdot))$ как $\Set*{\widetilde{a}^{-1/s} \widetilde{b}^{-1/s} f}{f \in L_s^P(u(\cdot, \cdot))}$ с нормой $\norm*{g}_{L_s^P(w)} = \norm*{g a^{1/s} b^{1/s}}_{L_s^P(u)}$, где, как и раньше, $\widetilde{a}, \widetilde{b}$ --- внешние функции с $\abs*{\widetilde{a}} = a, \abs*{\widetilde{b}} = b$. В чисто-решеточных терминах это выглядит так: мы взяли решетку $X = L^s$, добавили к ней вес $u^{-1/s}$, получив $X(u^{-1/s})$, после этого взяли подпространство, ``вырезанное'' проектором $P$, заданным на многочленах, получили $(X(u^{-1/s}))^P$, после чего добавили к получившейся решетке дополнительный вес $(\widetilde{a}\widetilde{b})^{-1/s}$, получив $(X(u^{-1/s}))^P((\widetilde{a}\widetilde{b})^{-1/s})$. \\

						Заметим также, что обозначения весов в дальнейшем могут отличаться от использованных в формулировках теорем выше. В основном, это сделано для того, чтобы читателю было легче сравнивать доказательства некоторых теорем-обобщений с их прообразами, выделять новые идеи. \\

						Наконец, хочется отметить, что тонкий вопрос определения пространств Харди здесь можно было бы вовсе обойти, потребовав, чтобы все веса были ограничены и отделены от нуля. В таком случае, они (пространства) определяются естественно и недвусмысленно, а результаты теорем имеют ценность, состоящую в том, что все полученные оценки будут зависеть только от $A_p$-констант и $BMO$-норм логарифмов весов, а не их (весов) существенных супремумов и инфимумов.
\newpage

					\section*{\begin{center}\S2 $K$-замкнутость при $0 < r \leq 1 < p < \infty$\end{center}}

\begin{theorem}
	\label{one_naib}
	Если веса $w_1(\cdot, \cdot)$, $w_2(\cdot, \cdot)$ удовлетворяют условиям
	\begin{enumerate}[label={\arabic*)}]
		\item существует такое $l \geq 1$, что $w_1(z_1, \cdot) \in A_l$ равномерно (т.е. ``$A_l$-константа'' ограничена как функция от $z_1$),
		\item $w_2(z_1, \cdot) \in A_p$ равномерно (в том же смысле),
		\item существует такое $m \geq 1$, что $w_1(\cdot, z_2), w_2(\cdot, z_2) \in A_m$ равномерно (т.е. ``$A_m$-константа'' ограничена как функция от $z_2$),
	\end{enumerate}
	то пара $\left(H^r(w_1(\cdot, \cdot)), H^p(w_2(\cdot, \cdot))\right)$ K-замкнута в $\left(L^r(w_1(\cdot, \cdot)), L^p(w_2(\cdot, \cdot))\right)$.
\end{theorem}
\begin{remark}
	Хватит, например, такого более обозримого (но более обременительного) условия:
	$$w_1(\cdot, \cdot) \in A_\infty,\qquad w_2(\cdot, \cdot) \in A_p,$$
	где, подчеркнем, фигурирующие выше условия $A_s$ --- двумерные. Это и составляет существо упомянутой во введении теоремы \ref{lft}.
\end{remark}
\begin{proof}
	Пусть $f = g + h$, где $f \in H^r(w_1(\cdot, \cdot)) + H^p(w_2(\cdot, \cdot))$, $g \in L^r(w_1(\cdot, \cdot))$, $h \in L^p(w_2(\cdot, \cdot))$. Так как, по условию, все веса удовлетворяют каким-то условиям Макенхапута равномерно по второй переменной, то для любого фиксированного значения $z_1 \in \mathbb{T}$ пространства $L^r(w_1(z_1, \cdot))$, $L^p(w_2(z_1, \cdot))$ являются $BMO$-регулярными квази-банаховыми решетками измеримых функций. Отсюда, вследствие общей теории (она изложена, например, в \cite{kis2}), мы можем сделать разложение функции $f$ аналитическим по второй переменной. Получим $f = g' + h'$, где $g' \in \ts$, $h' \in \ts$ и $\norm*{g'}_{L^r(w_1(\cdot, \cdot))} \lesssim \norm*{g}_{L^r(w_1(\cdot, \cdot))}$, $\norm*{h'}_{L^p(w_2(\cdot, \cdot))} \lesssim \norm*{h}_{L^p(w_2(\cdot, \cdot))}$. \\

	\textbf{Замечание}: Второй раз проделать тот же самый трюк нам мешает то, что полученные пространства $L^r(w_1(\cdot, \cdot))_A, L^p(w_2(\cdot, \cdot))_A$ не являются квази-банаховыми решетками измеримых функций, так как не удовлетворяет условию $f \in X, \abs*{g} \leq \abs*{f} \implies g \in X$. Индекс $A$ означает подпространство, состоящее из функций, аналитических по одной из переменных (в нашем случае по второй), это обозначение взято из \cite{kis2}. По поводу определения и свойств ``квази-банаховых решеток измеримых функций'' лучше всего см. \cite{GamKis}, хотя понятие является классическим и хорошо раскрыто в книге \cite{KantAk} (хотя оно там и не имеет нашего названия). Но все-таки эти пространства можно ``записать'' как квази-банаховы решетки измеримых функций, перейдя к коэффициентам разложения функций по общему безусловному базису, в этом-то и состоит главная идея доказательства. \\

	Теорема 3.10 из \cite{uncond}, благодаря условиям Макенхаупта, наложенным на веса по второй переменной, дает нам то, что для достаточно больших $k$, $k$-сплайновые вейвлеты одновременно образуют безусловный базис в пространствах $ReH^r(w_1(z_1, \cdot))$ и $ReH^p(w_2(z_1, \cdot))$ (с зафиксированной первой переменной). Нужно заметить, что в \cite{uncond} речь идет о вещественных пространствах Харди на $\mathbb{R}$. Переход от $ReH(\mathbb{R})$ к $ReH(\mathbb{T})$, как замечают сами авторы, осуществляется простой заменой обозначений. Таким образом, модифицированная для тора теорема из \cite{uncond} устанавливает наличие безусловного базиса (периодических) сплайновых вейвлетов в весовых $ReH^r$, $r\leq1$ и $L^p$, $p>1$, так как $ReH^p=L^p$ при $p>1$ (последнее соотношение верно для весовых пространств в силу теоремы 1 на странице 86 в \cite{Strom}). Уточним, пользуясь случаем, определение (ниже $\phi_r$ --- ядра Пуассона)
	$$ReH^s(w(\cdot)) = \Set*{f \text{--- комплекснозначное распределение на } \mathbb{T}}{f*\phi_r \in L^s(w(\cdot)),\ 0<r<1}$$
	с комплексными скалярами и нормой $\norm*{f}_{ReH^p(w(\cdot))} = \norm*{(\sup\limits_{0<r<1} |f*\phi_r|)(\cdot)}_{L^p(w(\cdot))}$. \\

	Нам понадобится определить некое вспомогательное пространство $ReH^s_2(w(\cdot,\cdot))$, состоящее из функций $f:\mathbb{T} \to D'$ ($D'$ обозначает распределения на торе), получающееся как замыкание по норме (в случае $s \geq 1$) или квази-норме (в случае $s < 1$) 
	$$\left(\int\norm*{f(z_1)}^s_{ReH^s(w(z_1, \cdot))}dz_1\right)^{\frac{1}{s}}$$
	множества конечных комбинаций вида $\sum\limits_{i=1}^n a_i(z_1) \chi_i(z_2)$, где $\{a_i(\cdot)\}$ -ограниченные измеримые функции. \\
	
	Докажем техническую лемму, которая позволит нам перемещаться между пространствами $L^s(w(\cdot, \cdot))_A$ и $ReH^s_2(w(\cdot,\cdot))$. \\

	\begin{lemma*}
		Если вес $w:\mathbb{T}\to\mathbb{R}_+$ удовлетворяет условию $A_\infty$, выполнено $r \leq 1$, то для аналитических функций на $\mathbb{T}$ (аналитичность понимается в смысле спектра) нормы $L^r(w)$ и $ReH^r(w)$ эквивалентны.
	\end{lemma*}
	\begin{remark}
		В реальности, нам нужна будет эквивалентность норм в пространствах $L^r(u(z_1, \cdot))$ и $ReH^r(u(z_1, \cdot))$, для которой константы эквивалентности не зависят от переменной $z_1$. В таком случае надо будет требовать $A_\infty$ в следующем смысле. Нужно, чтобы существовало такое число $\alpha < \infty$, что для любого $z_1$ вес $u(z_1, \cdot)$ удовлетворяет условию $A_\alpha$, причем $A_\alpha$-константу можно выбрать не зависящей от $z_1$.
	\end{remark}

	\begin{proof}
		Пусть $ReH^r_+(w)$ --- такое подпространство пространства $ReH^r(w)$, что у всех его элементов коэффициенты Фурье с отрицательными номерами равны нулю (то есть $ReH^r_+(w)$ это подпространство, состоящее из аналитических, в смысле спектра, распределений). \\

		Сворачивая распределения из пространства $ReH^r_+(w)$ с ядром Пуассона, будем получать функции в $H^r(w,\mathbb{D})$, пространстве функций, аналитических в круге, с ограниченными весовыми средними значениями по концентрическим окружностям. Такое преобразование является изоморфизмом. \\
		
		Теперь берем радиальные пределы получившихся функций в круге. Такое отображение также будет изоморфизмом. Действительно, норма граничной функции оценивается через норму в $H^r(w,\mathbb{D})$, благодаря лемме Фату. Докажем на плотном множестве, состоящем из ``хороших'' функций (можно взять на эту роль пересечение с $ReH^r_+(w)$ какого-нибудь из плотных в $ReH^r$ множеств, состоящих из гладких функций), обратную оценку
		$$\sup\limits_{0<\rho<1} \left(\int \abs*{(f*\phi_\rho)(z)}^r w(z)dz\right)^{\frac{1}{r}} \lesssim \norm*{f}_{L^r(w)}.$$
		Выберем натуральное число $m$ так, чтобы $r m > \alpha$, где $\alpha$ такое, что $w$ удовлетворяет условию $A_\alpha$ (такое $\alpha$ существует, так как вес удовлетворяет условию $A_\infty$). Заметим, что при таком выборе $m$, ядро Пуассона будет действовать в $L^{r m}(w)$. Профакторизуем функцию $f$ в произведение $b u$ внутренней и внешней функций (будет выполнено $\abs*{u} = \abs*{f}$). Тогда 
		\begin{multline*}
		\sup\limits_{0<\rho<1} \left(\int \abs*{(f*\phi_\rho)(z)}^r w(z)dz \right)^{\frac{1}{r}}
		\leq
		\sup\limits_{0<\rho<1} \left(\int \abs*{(u*\phi_\rho)(z)}^r w(z)dz \right)^{\frac{1}{r}}
		= \\
		=
		\sup\limits_{0<\rho<1} \left(\int \abs*{(u^{\frac{1}{m}}*\phi_\rho)(z)}^{rm} w(z)dz \right)^{\frac{1}{r}}
		\lesssim
		\norm*{u^{\frac{1}{m}}}_{L^{r m}(w)}^m
		=
		\norm*{u}_{L^{r}(w)}.
		=
		\norm*{f}_{L^{r}(w)}.
		\end{multline*}
		Таким образом, получили желанный изоморфизм. \\
		
		При этом, если мы начинали с настоящей функции, а не с распределения, то и придем к той же самой функции, а норма будет уже другая. Лемма доказана. \\
	
	\end{proof}

	Имеем, благодаря лемме, $g' \in ReH^r_2(w_1(\cdot,\cdot))$, $h' \in ReH^p_2(w_2(\cdot,\cdot))$, $\norm*{g'}_{ReH^r_2(w_1(\cdot,\cdot))} \lesssim \norm*{g'}_{L^r(w_1)}$, $\norm*{h'}_{ReH^p_2(w_2(\cdot,\cdot))} \lesssim \norm*{h'}_{L^p(w_2)}$. Причем соотношения для $g'$ не вытекают из леммы, а верны непосредственно, так как в этом случае пространство $ReH^p_2(w_2(\cdot,\cdot))$ ``вырождается'' в $L^p(w_2)$. \\

	Пусть $\left\{\chi_k(\cdot)\right\}_{k \in \mathbb{N}}$ --- тот самый базис (аргументом функций $\chi_k$ всегда будет $z_2$, заметим, что $\chi_k$ это не распределения, а функции). Теперь получим некую характеризацию нашего пространства $L^p(w_2(\cdot, \cdot))$.

	\begin{multline*}
		\norm*{\sum\limits_{i=1}^n a_i(z_1) \chi_i(z_2)}_{L^p(w_2(\cdot, \cdot))}
			\stackrel[(1)]{}{\sim}
		\frac{1}{2^n} \sum\limits_{\varepsilon_1, \dots, \varepsilon_n} \norm*{\sum\limits_{i=1}^n \varepsilon_i a_i(z_1) \chi_i(z_2)}_{L^p(w_2(\cdot, \cdot))}
			= \\
			=
		\int\limits_0^1 \norm*{\sum\limits_{i=1}^n r_i(t) a_i(z_1) \chi_i(z_2)}_{L^p(w_2(\cdot, \cdot))} dt
			\stackrel[(2)]{}{\sim}
		\left(\int\limits_0^1 \norm*{\sum\limits_{i=1}^n r_i(t) a_i(z_1) \chi_i(z_2)}^p_{L^p(w_2(\cdot, \cdot))} dt\right)^\frac{1}{p}
			= \\
			=
		\left(\int\limits_0^1 \int\limits_{\mathbb{T}}\int\limits_{\mathbb{T}} \abs*{\sum\limits_{i=1}^n r_i(t) a_i(z_1) \chi_i(z_2)}^p w_2(z_1, z_2) dz_2 dz_1 dt\right)^\frac{1}{p}
			= \\
			=
		\left(\int\limits_{\mathbb{T}}\int\limits_{\mathbb{T}} \int\limits_0^1 \abs*{\sum\limits_{i=1}^n r_i(t) a_i(z_1) \chi_i(z_2)}^p w_2(z_1, z_2) dt dz_2 dz_1 \right)^\frac{1}{p}
			\stackrel[(3)]{}{\sim} \\
			\stackrel[(3)]{}{\sim}
			\left(\int\limits_{\mathbb{T}}\int\limits_{\mathbb{T}} \left(\sum\limits_{i=1}^n \abs*{a_i(z_1) \chi_i(z_2)}^2\right)^\frac{p}{2} w_2(z_1, z_2) dz_2 dz_1 \right)^\frac{1}{p}.
	\end{multline*}
	Пояснения к пронумерованным переходам:
	\begin{enumerate}
		\item
			В силу безусловности базиса, которая означает ограниченность норм операторов $S_{N, \varepsilon}$ ($\varepsilon$ --- последовательность $1$ и $-1$, $\chi_k^*(x)$ --- коэффициент при $\chi_k$ в разложении $x$ по базису $\{\chi_i\}$),
			$$S_{N, \varepsilon}(x) := \sum_{k=1}^N \varepsilon_k \chi_k^*(x) \chi_k,$$
			константой, не зависящей от $\varepsilon$, которая (ограниченность), фактически и установлена в \cite{uncond}. Здесь ключевой является та самая ``равномерность'' условия Макенхаупта, которую мы требовали в условиях теоремы. Имея равномерность, можно увидеть, что нормы операторов по второй переменной $S_{N, \varepsilon}$, вообще говоря зависящие от переменной $z_1$, можно ограничить сверху константами, не зависящими от переменной $z_1$. Оценки снизу в данном случае получаются из оценок сверху, достаточно просто посмотреть на определение операторов $S_{N, \varepsilon}$.
		\item
			Неравенство Кахана (``векторное'' неравенство Хинчина, см. \cite{Woj} стр. 95). Отметим, что хотя неравенство Кахана обычно формулируется только для банаховых пространств, доказать его для квази-банаховых пространств не составляет труда.
		\item
			Скалярное неравенство Хинчина.
	\end{enumerate}

	Таким образом $ReH^p_2(w_2(\cdot,\cdot))$ изоморфно
	$$X_p := \Set*{\{a_k(\cdot)\}_{k=1}^\infty - \text{измеримая на $\mathbb{N} \times \mathbb{T}$}}{ \left(\sum\limits_{i = 1}^\infty \abs*{a_i(z_1) \chi_i(z_2)}^2\right)^{\frac{1}{2}} \in L^p(w_2(\cdot, \cdot))}.$$

	Проделаем нечто похожее для пространства $ReH^r_2(w_1(\cdot,\cdot))$. Как и в прошлый раз эквивалентность норм будем проверять только для конечных сумм вида $\sum\limits_{i=1}^n a_i(z_1) \chi_i(z_2)$, мы можем так делать потому, что они плотны в рассматриваемом пространстве. Здесь это будет играть существенную роль, так как, воспользовавшись теоремой 4 на стр. 87 из \cite{Strom} и гладкостью функций $\chi_k$, будем заменять обычную норму $\norm*{\sum\limits_{i=1}^n a_i(z_1) \chi_i(z_2)}_{ReH^r(w_1(z_1,\cdot))}$ на 
	$$\norm*{\sum\limits_{i=1}^n a_i(z_1) \chi_i(z_2)}_{L^r(w_1(z_1,\cdot))} + \qquad \quad \norm*{H\left(\sum\limits_{i=1}^n a_i(z_1) \chi_i(z_2)\right)}_{L^r(w_1(z_1,\cdot))}.$$
	Приступим к выкладкам.

	\begin{multline*}
		\norm*{\sum\limits_{i=1}^n a_i(z_1) \chi_i(z_2)}_{ReH^r_2(w_1(\cdot,\cdot))}^r
			\sim
		\frac{1}{2^n} \sum\limits_{\varepsilon_1, \dots, \varepsilon_n} \norm*{\sum\limits_{i=1}^n \varepsilon_i a_i(z_1) \chi_i(z_2)}_{ReH^r_2(w_1(\cdot,\cdot))}^r
			= \\
			=
			\int\int \left( \int\limits_0^1 \abs*{\sum\limits_{i=1}^n r_i(t) a_i(z_1) \chi_i(z_2)}^r + \int\limits_0^1 \abs*{\sum\limits_{i=1}^n r_i(t) a_i(z_1) (H\chi_i)(z_2)}^r \right) w_1(z_1,z_2) dz_2 dz_1
			\sim \\
			\sim
			\int\int \left( \left(\sum\limits_{i=1}^n \abs*{a_i(z_1) \chi_i(z_2)}^2 \right)^{\frac{r}{2}} + \left(\sum\limits_{i=1}^n \abs*{a_i(z_1) (H\chi_i)(z_2)}^2 \right)^{\frac{r}{2}} \right) w_1(z_1,z_2) dz_2 dz_1.
	\end{multline*}

	Таким образом $ReH^r_2(w_1(\cdot,\cdot))$ изоморфно (ниже $\{a_k(\cdot)\}_{k=1}^\infty$ --- измеримые на $\mathbb{N} \times \mathbb{T}$)
	$$X_r := \Set*{\{a_k(\cdot)\}_{k=1}^\infty}{\left(\sum\limits_{i=1}^n \abs*{a_i(z_1) \chi_i(z_2)}^2 \right)^{\frac{1}{2}} + \left(\sum\limits_{i=1}^n \abs*{a_i(z_1) (H\chi_i)(z_2)}^2 \right)^{\frac{1}{2}} \in L^r(w_1(\cdot, \cdot))}.$$

	Несложно увидеть, что оба получившихся пространства будут квази-банаховыми решетками измеримых функций. \\
	
	Проверим, что решетка $X_p$ является $BMO$-регулярной. Для этого докажем, что оператор гармонического сопряжения $H$ действует в $(X_p)^\alpha$ ($\alpha$-конвексификации решетки $X_p$, определение можно найти в \cite{kis2}). \\

	Оценка, которую мы пытаемся получить, выглядит так: 
	$$\int\limits_{\mathbb{T}} \int\limits_{\mathbb{T}} \left(\sum\limits_{i=1}^\infty \abs*{(H a_i)(z_1)}^{2\alpha} \abs*{\chi_i(z_2)}^2\right)^\frac{p}{2} w_2(z_1, z_2) dz_1 dz_2
	\lesssim
	\int\limits_{\mathbb{T}} \int\limits_{\mathbb{T}} \left(\sum\limits_{i=1}^\infty \abs*{a_i(z_1)}^{2\alpha} \abs*{\chi_i(z_2)}^2\right)^\frac{p}{2} w_2(z_1, z_2) dz_1 dz_2.$$
	Ее несложно получить, воспользовавшись условием 3. Действительно,

	\begin{multline*}
		\int\limits_{\mathbb{T}} \left(\sum\limits_{i=1}^\infty \abs*{(H a_i)(z_1)}^{2\alpha} \abs*{\chi_i(z_2)}^2\right)^\frac{p}{2} w_2(z_1, z_2) dz_1 = \int\limits_{\mathbb{T}} \left(\sum\limits_{i=1}^\infty \abs*{(H (a_i \abs*{\chi_i(z_2)}^\frac{1}{\alpha}))(z_1)}^{2\alpha}\right)^\frac{p}{2} w_2(z_1, z_2) dz_1 \leq \dots
	\end{multline*}

	Как известно, $H$ --- оператор Кальдерона-Зигмунда на пространстве $L^{p\alpha}(l^{2\alpha})$, так что, в силу условия 3 (нужно выбрать $\alpha$ так, чтобы $p \alpha$ было больше $m$), можем продолжить выкладку выше
	$$\dots \leq C(z_1) \int\limits_{\mathbb{T}} \left(\sum\limits_{i=1}^\infty \abs*{a_i(z_1) \abs*{\chi_i(z_2)}^\frac{1}{\alpha}}^{2\alpha}\right)^\frac{p}{2} w_2(z_1, z_2) dz_1 \leq C \int\limits_{\mathbb{T}} \left(\sum\limits_{i=1}^\infty \abs*{a_i(z_1)}^{2\alpha} \abs*{\chi_i(z_2)}^2\right)^\frac{p}{2} w_2(z_1, z_2)dz_1.$$

	Выше мы также, пользуясь равномерностью наложенного условия Макенхаупта, заменили константу $C(z_1)$ на независимую от $z_1$ константу $C$.
	Осталось лишь проинтегрировать полученную оценку по переменной $z_2$. $BMO$-регулярность доказана. \\
	
	$BMO$-регулярность второй решетки доказывается совершенно аналогично. \\
	
	Возвратимся к нашим разложениям. Снова используя общую теорию, можем сделать коэффициенты $a_i(\cdot)$ аналитическими функциями. Таким образом получим разложение $f=g''+h''$, где $g''$ и $h''$ лежат в соответствующих $ReH^s_2(w_i(\cdot,\cdot))$, имеют аналитические коэффициенты в разложении по базису $\{\chi_i(\cdot)\}$, и их $ReH^s_2(w_i(\cdot,\cdot))$ нормы нужным образом оцениваются через $ReH^s_2(w_i(\cdot,\cdot))$ нормы $g'$ и $h'$. Заметим, что коэффициенты $g''$ и $h''$ были получены из $g'$, $h'$ домножением на некоторую функцию переменной $z_1$ (так работает ``общая теория''), так что $g''$ и $h''$ так и остались аналитическими по второй переменной. Использовав еще раз лемму, заканчиваем доказательство.
\end{proof}

\newpage
					\section*{\begin{center}\S3 $K$-замкнутость при $r > 1$, $p = \infty$\end{center}}
Как известно, $K$-замкнутость пары
$$\left(H_p(b_1(z_1) w_1(z_1, z_2) a_1(z_2)), H_\infty(b_2(z_1) w_2(z_1,z_2) a_2(z_2))\right)$$
в паре
$$\left(L_p(b_1(z_1) w_1(z_1, z_2) a_1(z_2)), L_\infty(b_2(z_1) w_2(z_1,z_2) a_2(z_2))\right)$$
эквивалентна $K$-замкнутости соответствующей пары ``преданнуляторов'' в соответствующей паре предсопряженных пространств. То есть $K$-замкнутости пары
$$\left(L_1^P(\widetilde{b}_1(z_1) \widetilde{w}_1(z_1, z_2) \widetilde{a}_1(z_2)), L_q^P(\widetilde{b}_2(z_1) \widetilde{w}_2(z_1,z_2) \widetilde{a}_2(z_2))\right)$$
в паре
$$\left(L_1(\widetilde{b}_1(z_1) \widetilde{w}_1(z_1, z_2) \widetilde{a}_1(z_2)), L_q(\widetilde{b}_2(z_1) \widetilde{w}_2(z_1,z_2) \widetilde{a}_2(z_2))\right),$$
где пространства $L_s^P(w)$ были определены во введении, верно соотношение $\frac{1}{p}+\frac{1}{q}=1$, а веса пересчитываются следующим образом:
$$\widetilde{b}_1(z_1) = b_2(z_1), \qquad \widetilde{w}_1(z_1, z_2) = w_2(z_1, z_2), \qquad \widetilde{a}_1(z_2) = a_2(z_2);$$
$$\widetilde{b}_2(z_1) = b^{1-q}_1(z_1), \qquad \widetilde{w}_2(z_1, z_2) = w^{1-q}_1(z_1, z_2), \qquad \widetilde{a}_2(z_2) = a^{1-q}_1(z_2).$$
В дальнейшем мы из соображений удобства не будем писать тильду над $b_1,w_1,a_1,b_2,w_2,a_2$, когда имеем дело с ``преданнуляторами'' и предсопряженными пространствами.

						\subsection*{\begin{center} Переход от двух весов к одному весу и окаймленному оператору\end{center}}
Для начала заметим, что без ограничения общности можно считать, что $a_1 = a_2 = a$, $b_1 = 1$, $b_2 = b$. Этого всегда можно добиться в самом начале, используя в исходной задаче домножение на нужные внешние функции. \\

\noindent Теперь, используя прием из \cite{ruts} (стр 192), получаем, что вопрос K-замкнутости пары
$$\left(L_1^P(w_1(z_1, z_2) a(z_2)), L_q^P(b(z_1) w_2(z_1,z_2) a(z_2))\right)$$
в паре
$$\left(L_1(w_1(z_1, z_2) a(z_2)), L_q(b(z_1) w_2(z_1,z_2) a(z_2))\right)$$
эквивалентен вопросу о K-замкнутости пары
$$\left(L_1^{P^u}(w(z_1, z_2) a(z_2)), L_q^{P^u}(b(z_1) w(z_1, z_2) a(z_2))\right)$$
в паре
$$\left(L_1(w(z_1, z_2) a(z_2)), L_q(b(z_1) w(z_1, z_2) a(z_2))\right),$$
при $w = \frac{w_1^{\frac{q}{q-1}}}{w_2^{\frac{1}{q-1}}}$, $u = \frac{w_1^{\frac{1}{q-1}}}{w_2^{\frac{1}{q-1}}}$, $P^u f = u^{-1}P(uf)$. Тут нужно отметить, что используемые здесь $L^{P_u}(w)$ определяются ``чисто решеточным'' образом, так что нам нет необходимости накладывать какие-то дополнительные условия на $w$ и $u$, чтобы такая запись имела смысл. Уточним, что если вспомнить ``решеточный процесс'', стоящий за определением пространств $L_s^P(w)$, то переход к окаймленному оператору будет совершен на втором шаге, перед добавлением ``разделяющегося'' веса, можно считать, что за окаймленным оператором закреплено плотное множество $u^{-1} \cdot D$, где $D$ --- множество тригонометрических многочленов.

						\subsection*{\begin{center}Основная теорема\end{center}}

\begin{theorem}
	K-замкнутость $\left(L_1^{P^u}(w(z_1, z_2) a(z_2)), L_q^{P^u}(b(z_1) w(z_1, z_2) a(z_2))\right)$ имеет место, если
	\begin{enumerate}[label={\arabic*)}]
		\item[0)] $w_1(\cdot, \cdot), w_2(\cdot, \cdot) \in L^1(\mathbb{T}^2)$ (техническое требование для корректности наших определений),
		\item  $w = \frac{w_1^{\frac{q}{q-1}}}{w_2^{\frac{1}{q-1}}}$ удовлетворяет условию $A_\infty$ по второй переменной равномерно (под равномерностью мы здесь понимаем существование константы в обратном неравенстве Гельдера, не зависящей от первой переменной),
		\item $w_1$ удовлетворяет условию $A_1$ по второй переменной равномерно,
		\item $\log(a), \log(b) \in \text{BMO}$,
		\item $\log(w(\cdot, z_2))$ лежит в пространстве $\text{BMO}$ по первой переменной равномерно.
		\item $w_2 \in A_q$ равномерно по второй переменной,
	\end{enumerate}
\end{theorem}
\begin{remark}
	Кажется, что при вышеуказанных условиях вес $b$ является избыточным, что его можно считать частью веса $w$ и все условия будут удовлетворены. Это иллюзия, ведь при определении пространств Харди двух переменных мы сразу же наложили условие, заставляющее все ``неразделяющиеся'' веса двух переменных лежать в $L^1(\mathbb{T}^2)$.
\end{remark}
\begin{proof}
	Введем обозначения: $Y_1 := L_1^{P^u}(w(z_1, z_2) a(z_2))$, $Y_2 := L_q^{P^u}(b(z_1) w(z_1, z_2) a(z_2))$. \\

	Пусть $Y_1 + Y_2 \ni f = g + h$, где $g \in L_1(w(z_1, z_2) a(z_2)), h \in L_q(b(z_1) w(z_1, z_2) a(z_2))$. Обозначим $\norm*{g}_{L_1(w(z_1, z_2) a(z_2))} =: A$, $\norm*{h}_{L_q(b(z_1) w(z_1, z_2) a(z_2))} =: B$ и будем искать такие функции $g' \in Y_1$, $h' \in Y_2$, что $f = g' + h'$ и выполнено $\norm*{g'}_{L_1(w(z_1, z_2) a(z_2))} \lesssim A$, $\norm*{h'}_{L^q(b(z_1) w(z_1, z_2) a(z_2))} \lesssim B$. \\

	Благодаря условию $\log a \in \text{BMO}$, можем найти такой набор (``аналитическое разложение единицы, подчиненное весу $a$'') $\{\phi_j\}_{j \in \mathbb{Z}} \subseteq H^\infty(\mathbb{T})$, что выполнено
	$$\abs*{\phi_j}^{\frac{1}{8}} a \lesssim 2^j, \sum \abs{\phi_j}^{\frac{1}{8}} 2^j \lesssim a$$
	$$\sum \abs*{\phi_j}^{\frac{1}{8}} \leq c, \sum \phi_j = 1.$$
	Построение см. в \cite{KisBurgProj}. \\

	Определим функции $\psi_j$ из факторизационного соотношения $\phi_j = \theta_j \psi_j^8$, где $\theta_j$ --- внутренняя функция, а $\psi_j$ --- внешняя. Имеем
	$$f = \sum \theta_j \psi_j^4 f \psi_j^4 = \sum \theta_j \psi_j^4 g \psi_j^4 + \sum \theta_j \psi_j^4 h \psi_j^4,$$
	а также $f \psi_j^4 = g \psi_j^4 + h \psi_j^4$, именно с этим разложением мы будем работать в ближайшее время. \\
	
	Пусть $P_2$ --- естественный проектор на $\bs$. Можно считать, что это оператор, действующий на функцию $v(z_1, \cdot)$ как $(I-\mathfrak{R})$, где $\mathfrak{R}$ --- это классический одномерный проектор Рисса (вопросов с измеримостью по двум переменным здесь не возникает). Уточним, что $P_2$ аннулирует все константы. Благодаря условиям 1 и 2, для любого фиксированного $z_1$ можем построить весовое разложение Кальдерона-Зигмунда для функции $g \psi_j^4$ по второй переменной для оператора $P_2^u$ по уровню $\lambda(z_1)$, который мы выберем позже (точную формулировку используемой теоремы см. в \cite{ruts} (стр. 191), а ее доказательство см. в \cite{KisAnis} (стр. 761-763)), т.е. существуют такие функции $g_0^j(z_1, \cdot)$, $g_1^j(z_1, \cdot)$, и такое семейство множеств $\{\Omega_{z_1}^j\}_{z_1 \in \mathbb{T}}$, что:
	\begin{itemize}
		\item
			$(g\psi_j^4)(z_1, \cdot) = g_0^j(z_1, \cdot) + g_1^j(z_1, \cdot)$,
		\item
			$\abs*{g_0^j(z_1, \cdot)} \lesssim \lambda_j(z_1)$,
		\item
			$\int \abs*{g_0^j(z_1, z_2)} w(z_1, z_2) dz_2 \lesssim \int \abs*{g \psi_j^4(z_1, z_2)} w(z_1, z_2) dz_2$
		\item
			$\int \abs*{g_1^j(z_1, z_2)} w(z_1, z_2) dz_2 \lesssim \int \abs*{g \psi_j^4(z_1, z_2)} w(z_1, z_2) dz_2$
		\item
			$\text{supp}(g_1^j(z_1, \cdot)) \subseteq \Omega_{z_1}^j$, где $(w(z_1, \cdot) \mu)(\Omega_{z_1}^j) \lesssim \frac{\int \abs*{g \psi_j^4(z_1, z_2)} w(z_1, z_2) dz_2}{\lambda_j(z_1)}$,
		\item
			$\int\limits_{\mathbb{T} \setminus \Omega_{z_1}^j} \abs{(P_2^u g_1^j)(z_1, z_2)} w(z_1, z_2) dz_2 \leq \int \abs*{g \psi_j^4(z_1, z_2)} w(z_1, z_2) dz_2$.
	\end{itemize}

	Здесь нужно заметить, что все константы, скрытые в оценках выше за символом ``$\lesssim$'', не зависят от переменной $z_1$. Это достигается благодаря требованию равномерности условий Макенхаупта для весов $w$ и $w_1$ и специфики теоремы 4 из \cite{KisAnis}. \\

	Займемся теперь определением чисел $\lambda_j(z_1)$. Пусть $y_j := \left(\int \abs*{(h \psi_j^4)(z_1, z_2)}^q w(z_1, z_2) dz_2\right)^{\frac{1}{q}}$. Тогда, используя свойства аналитического разложения единицы, получим
	\begin{gather*}
	\sum\limits_j 2^j y_j(z_1)^q = \sum\limits_j 2^j \int \abs*{(h \psi_j^4)(z_1, z_2)}^q w(z_1, z_2) dz_2 \lesssim \int \abs*{h(z_1, z_2)}^q \left(\sum\limits_j 2^j \abs*{\psi_j} \right) w(z_1, z_2) dz_2 \lesssim \\
	\lesssim \int \abs*{h(z_1, z_2)}^q w(z_1, z_2) a(z_2) dz_2,
	\end{gather*}
	откуда сразу видно, что
	$$\left( \int \sum\limits_j 2^j y_j(z_1)^q b(z_1) dz_1 \right)^{\frac{1}{q}} \lesssim B.$$

	В левой части неравенства выше стоит норма последовательности $\{2^{\frac{j}{q}}y_j(\cdot)\}_{j \in \mathbb{Z}}$ измеримых функций на $\mathbb{T}$ в решетке $L^q(l^q, b)$, которая (решетка) $BMO$-регулярна, так как $\log(b) \in \text{BMO}$, так что мы можем найти такие функции $v_j \geq y_j$, что $\sup\limits_j \norm*{\log v_j}_{BMO} \leq c$ и $\left( \int \sum\limits_j 2^j v_j(z_1)^q b(z_1) dz_1 \right)^{\frac{1}{q}} \lesssim B$. Так как $\log(w(\cdot, z_2)) \in \text{BMO}$ равномерно (по условию 4) и $\sup\limits_j \norm*{\log v_j}_{BMO} \leq c$, то, по лемме 7.2 из \cite{GamKis}, можем найти такие функции $\gimel_j(z_1, z_2)$, что, во-первых, $v_j(z_1) w(z_1, z_2) \lesssim \gimel_j(z_1, z_2) \lesssim v_j(z_1) w(z_1, z_2)$, и, во-вторых, $\abs*{H(\gimel_j^{\frac{1}{k}})} \lesssim \gimel_j^{\frac{1}{k}}$ для какого-то числа $k \geq 2$ (с абсолютными константами внутри символа ``$\lesssim$'' и показателем $k$, не зависящем от индекса $j$, все благодаря равномерности). \\
	
	Теперь, наконец,
	$$\lambda_j(z_1) := \frac{v_j(z_1)^p}{\left( \int \abs*{(g \psi_j^4)(z_1, z_2)} w(z_1, z_2) dz_2\right)^{p-1}}.$$

	Сразу же
	\begin{gather*}
	\int \abs*{g_0^j (z_1, z_2)}^q w(z_1, z_2) dz_2 \lesssim \lambda_j(z_1)^{q-1} \int \abs*{(g \psi_j^4)(z_1, z_2)} w(z_1, z_2) dz_2 = \\
	= v_j(z_1)^{p(q-1)}  \left( \int \abs*{(g \psi_j^4)(z_1, z_2)} w(z_1, z_2) dz_2 \right)^{1-(p-1)(q-1)} = v_j(z_1)^q,
	\end{gather*}
	то есть
	$$\left( \int \abs*{g_0^j (z_1, z_2)}^q w(z_1, z_2) dz_2 \right)^{\frac{1}{q}} \lesssim v_j(z_1).$$

	Теперь положим (используемые ниже функции $\Phi_j$, аналитичные по первой переменной, мы выберем позже)
	$$u_j := P_2^u(f \psi_j^4) = P_2^u(g_1^j) + P_2^u(g_0^j + h \psi_j^4),$$
	$$\alpha_j := \Phi_j u_j - P_2^u(g_0^j + h \psi_j^4),$$
	$$\Lambda := \sum\limits_j \theta_j \psi_j^4 \alpha_j.$$

	Тогда $f = \sum\limits_j \theta_j \psi_j^4 f \psi_j^4 = (\sum\limits_j \theta_j \psi_j^4 g_1^j - \Lambda) + (\sum\limits_j \theta_j \psi_j^4 (g_0^j + h \psi_j^4) + \Lambda) = g' + h'$, где $g'$ и $h'$ --- искомые. Докажем это. \\

	Проверим, что $P^u(h')=h'$ (утверждение $P^u(g')=g'$ сразу же последует из этого). Во-первых, легко убедиться в том, что
	$$h' = (I-P_2^u)(g_0^j+h \psi_j^4) + \Phi_j P_2^u(f \psi_j^4).$$
	Выше мы сразу исключили из $h'$ сумму и множители $\theta_j \psi_j^4$, легко понять, что это не ограничивает общности. Более того, дальше мы будем считать, что $h'$ оказалась функцией вида $u^{-1} v$, где $v$ --- тригонометрический многочлен, это предположение нас также нисколько не ограничивает, ведь распространение на общий случай легко получается из соображений плотности таких функций. Теперь просто формальным образом проверим принадлежность пространствам $L^{P^u}$ для каждого слагаемого. Для проверки первого слагаемого достаточно провести прямое вычисление:
	\begin{gather*}
	P^u(I-P_2^u)(g_0^j+h\psi_j^4)
	=
	u^{-1} P u u^{-1} (I-P_2) u (g_0^j+h\psi_j^4)
	=
	u^{-1} P (I-P_2) u (g_0^j+h\psi_j^4)
	= \\ =
	u^{-1} (I-P_2) u (g_0^j+h\psi_j^4)
	=
	(I-P_2^u)(g_0^j+h\psi_j^4).
	\end{gather*}
	Для проверки второго будем действовать немного иначе. Знаем, что $P^u (f) = f$, т.е. $P(uf) = uf$, т.е. $uf \in \nlbs$, откуда очевидно, что $uf \psi_j \in \nlbs$, следовательно $P_2 (u f \psi_j) \in \rbs$, значит $\Phi_j P_2 (u f \psi_j) \in \rs$, откуда, наконец, $P^u \Phi_j P_2^u (f \psi_j) = u^{-1} P (\Phi_j P_2 (u f \psi_j)) = u^{-1} \Phi_j P_2 (u f \psi_j) = \Phi_j P_2^u (f \psi_j)$, что и требовалось доказать. \\

	Теперь нужно проверить оценки. Для начала установим нужные оценки для членов, не содержащих $\Lambda$:
	\begin{gather*}
		\left( \int \int \abs*{\sum \theta_j \psi_j^4 (g_0^j + h \psi_j^4)}^q b(z_1) w(z_1, z_2) a(z_2) dz_2  dz_1 \right)^{\frac{1}{q}} \lesssim \\
		\lesssim \left( \int \left(\sum \int  2^j \abs*{g_0^j + h \psi_j^4}^q w(z_1, z_2) dz_2 \right) b(z_1) dz_1 \right)^{\frac{1}{q}} \lesssim \tag{$\star$} \\
		\lesssim \left(\int \left(\sum 2^j v_j^q + \int \left( \abs*{h}^q \sum 2^j \abs*{\psi_j} w(z_1, z_2) dz_2 \right) \right) b(z_1) dz_1\right)^{\frac{1}{q}} \lesssim \\
		\lesssim \left( \int \sum 2^j v_j^q b(z_1) dz_1 \right)^{\frac{1}{q}} + \left( \int \int \abs*{h}^q b(z_1) w(z_1, z_2) a(z_2) dz_2 dz_1 \right)^{\frac{1}{q}} \lesssim B. \tag{$\star\star$}
	\end{gather*}
	Аналогичным образом, но еще проще:
	\begin{gather*}
		\int \int \abs*{\sum \theta_j \psi_j^4 g_1^j}  w(z_1, z_2) a(z_2) dz_2 dz_1 \lesssim
		\int \int \sum \abs*{g_1^j} \abs*{\psi_j} a(z_2) w(z_1, z_2) dz_2  dz_1 \lesssim \\
		\lesssim \int \left( \sum 2^j \left( \int \abs*{g_1^j} w(z_1, z_2) dz_2 \right) \right)  dz_1 \lesssim
		\int \left( \int \abs*{g} \sum 2^j \abs*{\psi_j} w(z_1, z_2) dz_2 \right)  dz_1 \lesssim \\
		\lesssim \int \int \abs*{g}  w(z_1, z_2) a(z_2) dz_2 dz_1 = A.
	\end{gather*}
	
	Переходим к выбору функций $\Phi_j$. Для начала, введем обозначение $r_j(z_1, z_2) := \gimel_j^\frac{1}{k}$, а также найдем такое число $s \in \mathbb{N}$, что $p/s < 1$. Теперь:
	$$\gamma_j(z_1, z_2) := \max \left\{ 1, \left(\frac{\abs*{(P_2^u g_1^j) (z_1, z_2)}}{\lambda_j(z_1)}\right)^{\frac{1}{k s}} \right\}$$
	$$F_j := \frac{r_j + i H(r_j)}{r_j \gamma_j + i H(r_j \gamma_j)},\qquad \qquad \Phi_j = 1 - (1 - F_j^{k s})^k .$$
	Тогда $\Phi_j$ --- аналитичны по переменной $z_1$ ($H$ здесь действует по $z_1$). Также выполнена оценка $\abs*{\Phi_j} \lesssim \frac{1}{\gamma_j^{k s}}$, так как 
	$$\abs*{F_j} = \abs*{ \frac{r_j + i H(r_j)}{r_j \gamma_j + i H(r_j \gamma_j)}} \leq \abs*{ \frac{r_j + i H(r_j)}{r_j \gamma_j}} \lesssim \frac{\abs*{r_j} + \abs*{r_j}}{r_j \gamma_j} \lesssim \frac{1}{\gamma_j} \leq 1,$$
	$$\abs*{\Phi_j} = \abs*{- C_k^1 F_j^{ks} + C_k^2 F_j^{2ks} - \dots \pm C_k^k F_j^{k^2 s}} \stackrel[\text{т.к. $\abs*{F_j} \lesssim 1$}]{}{\lesssim} \abs*{F_j^{ks}} \lesssim \frac{1}{\gamma_j^{k s}},$$
	а также имеет место $\abs*{\Phi_j}\abs*{P_2^u(g_1^j)} \lesssim \frac{\abs*{P_2^u(g_1^j)}}{\gamma_j^{k s}} \leq \lambda_j(z_1)$. Все константы, скрытые здесь под символом ``$\lesssim$'', абсолютные в смысле независимости от индекса $j$. \\

	Проверим неравенство $\left( \int \int \abs*{\Lambda}^q b(z_1) w(z_1, z_2) a(z_2) dz_2 dz_1 \right)^{\frac{1}{q}} \lesssim B$. В силу оценки $\Phi_j \lesssim 1$, запишем
	\begin{gather*}
		\left( \int \int \abs*{\Lambda}^q b(z_1) w(z_1, z_2) a(z_2) dz_2 dz_1 \right)^{\frac{1}{q}} = \\
		= \left( \int \int \abs*{\sum \theta_j \psi_j^4 \left( \Phi_j P_2^u(g_1^j) + \Phi_j P_2^u(g_0^j+ h \psi_j^4) - P_2^u(g_0^j + h \psi_j^4) \right) }^q b(z_1) w(z_1, z_2) a(z_2) dz_2 dz_1 \right)^{\frac{1}{q}} \lesssim \\
		\lesssim \left( \int \int \abs*{\sum \theta_j \psi_j^4 \Phi_j P_2^u(g_1^j)}^q b(z_1) w(z_1, z_2) a(z_2) dz_2 dz_1 \right)^{\frac{1}{q}} + \\
		+ \left( \int \int \abs*{\sum \theta_j \psi_j^4 P_2^u(g_0^j + h \psi_j^4)}^q b(z_1) w(z_1, z_2) a(z_2) dz_2 dz_1 \right)^{\frac{1}{q}} \lesssim \\
		\lesssim \left( \int \int \sum 2^j \abs*{\Phi_j P_2^u(g_1^j)}^q b(z_1) w(z_1, z_2) dz_2 dz_1 \right)^{\frac{1}{q}} + \\
		+ \left( \int \int \sum 2^j \abs*{P_2^u(g_0^j + h \psi_j^4)}^q b(z_1) w(z_1, z_2) dz_2 dz_1 \right)^{\frac{1}{q}} = I_1 + I_2
	\end{gather*}

	Чтобы оценить интеграл $I_2$, достаточно воспользоваться сильным типом $(q,q)$ оператора $P_2^u$ относительно веса $w$ (наличие сильного типа $(q,q)$ очевидно следует из условия 5), а далее провести рассуждения полностью совпадающие с ``хвостом'' (см. $\star$ -- $\star\star$) оценки, проделанной выше. \\

	Чтобы оценить интеграл $I_1$, сначала оценим внутренний интеграл $\int \abs*{\Phi_j P_2^u g_1^j}^q w(z_1, z_2) dz_2$ для фиксированного значения $z_1$. Пусть
	$$\rho_{j, z_1}(t) = (w(z_1, \cdot) \mu) \Set*{z_2}{\abs*{\Phi_j P_2^u g_1^j} > t},$$
	тогда, в силу оценки $\abs*{\Phi_j}\abs*{P_2^u(g_1^j)} \leq c \lambda_j(z_1)$ и слабого типа $(1,1)$ оператора $P_2^u$ относительно веса $w$, который получается из теоремы 4 из \cite{KisAnis}, имеем
	\begin{gather*}
	\int \abs*{\Phi_j P_2^u g_1^j}^q w(z_1, z_2) dz_2 \lesssim
	\int\limits_0^{c \lambda_j(z_1)} t^{q-1} \rho_{j,z_1}(t)dt \lesssim
	\lambda_j(z_1)^{q-1} \int \abs*{g \psi_j^4} w(z_1, z_2) dz_2 = \\
	= \left( \lambda_j(z_1) \left( \int \abs*{g \psi_j^4} w(z_1, z_2) dz_2 \right)^{p-1} \right)^{\frac{1}{p-1}} \leq v_j(z_1)^{\frac{p}{p-1}} = v_j(z_1)^q
	\end{gather*}

	Отсюда видно, что
	$$I_1 \lesssim \left( \int \sum 2^j v_j^q b(z_1) dz_1 \right)^{\frac{1}{q}} \lesssim B.$$

	Осталось проверить последнее неравенство $\int \int \abs*{\Lambda} w(z_1, z_2) a(z_2) dz_2 dz_1 \lesssim A$. Легко видеть (подобное мы уже делали выше), что
	\begin{gather*}
		\int \int \abs*{\Lambda} w(z_1, z_2) a(z_2) dz_2 dz_1 \lesssim 
		\int \left( \sum 2^j \int \abs*{\Phi_j P_2^u(g_1^j)} w(z_1, z_2) dz_2 dz_1 \right) + \\
		+ \int \left( \sum 2^j \int \abs*{1-\Phi_j}\abs*{P_2^u(g_0^j + h \psi_j^4)} w(z_1, z_2) dz_2 dz_1 \right) =
		I_3 + I_4.
	\end{gather*}

	Оценим интеграл $I_3$. Для начала,
	$$\int \abs*{\Phi_j P_2^u g_1^j} w(z_1, z_2) dz_2 \lesssim \lambda_j(z_1) (w(z_1, \cdot)\mu)(\Omega_{z_1}^j) + \int\limits_{\mathbb{T} \setminus \Omega_{z_1}^j} \abs{P_2^u g_1^j} w(z_1, z_2) dz_2 \lesssim \int \abs*{g \psi_j^4} w(z_1, z_2) dz_2,$$
	откуда уже хорошо знакомым приемом получаем $I_3 \lesssim A$.

	Вспомним неравенства $\left(\int \abs*{(h \psi_j^4)(z_1, z_2)}^q w(z_1, z_2) dz_2\right)^{\frac{1}{q}} \leq v_j(z_1)$ и $\left( \int \abs*{g_0^j (z_1, z_2)}^q w(z_1, z_2) dz_2 \right)^{\frac{1}{q}} \lesssim v_j(z_1)$, а также воспользуемся сильным типом $(q,q)$ у оператора $P_2^u$ относительно веса $w(z_1, \cdot)$, получим
	\begin{gather*}
	I_4 \lesssim
	\int \left( \sum 2^j \left( \int \abs*{1-\Phi_j}^p w(z_1, z_2) dz_2 \right)^{\frac{1}{p}} \left( \int \abs*{P_2^u(g_0^j + h \psi_j^4)}^q w(z_1, z_2) dz_2 \right)^{\frac{1}{q}} dz_1 \right)
	\lesssim \\ \lesssim
	\int \left( \sum 2^j \left( \int \abs*{1-\Phi_j}^p w(z_1, z_2) dz_2 \right)^{\frac{1}{p}} \left( \int \abs*{g_0^j + h \psi_j^4}^q w(z_1, z_2) dz_2 \right)^{\frac{1}{q}} dz_1 \right)
	\lesssim \\ \lesssim
	\int \left( \sum 2^j \left( \int \abs*{1-\Phi_j}^p w(z_1, z_2) dz_2 \right)^{\frac{1}{p}} v_j dz_1 \right)
	\end{gather*}

	Очевидно, что выполнено
	$$\abs*{1-\Phi_j} = \abs*{1-F_j^{k s}}^k = \abs*{1-F_j}^k \abs*{1+F_j+F_j^2+\dots+F_j^{ks-1}}^k \lesssim \abs*{1-F_j}^k,$$
	$$\abs*{1-F_j} = \abs*{\frac{r_j (\gamma_j - 1) + i H(r_j (\gamma_j - 1))}{r_j \gamma_j + i H(r_j \gamma_j)}} \lesssim \frac{\gamma_j - 1}{\gamma_j} + \frac{\abs*{H(r_j(\gamma_j-1))}}{r_j},$$
	последнее выражение равняется $\frac{\abs*{H(r_j(\gamma_j-1))}}{r_j}$ на множестве, где $\gamma_j = 1$, т.е. на множестве, где $\abs*{P_2^u g_1^j} \leq \lambda_j(z_1)$. Отсюда
	\begin{gather*}
	I_4 \lesssim \int \sum 2^j \left((w(z_1, \cdot) \mu)\left(\Set*{z_2}{\gamma_j(z_1, z_2) > 1}\right)\right)^{\frac{1}{p}} v_j(z_1) dz_1
	+ \\ +
	\int \sum 2^j \left( \int \abs*{\frac{H(r_j (\gamma_j-1))^{k p}}{r_j^{k p}}} w(z_1, z_2) dz_2 \right)^{\frac{1}{p}} v_j(z_1) dz_1.
	\end{gather*}
	Вспоминая, что $r_j^k = \gimel_j$, $v_j(z_1) w(z_1, z_2)^{\frac{1}{p}} \lesssim \gimel_j(z_1, z_2) \lesssim v_j(z_1) w(z_1, z_2)^{\frac{1}{p}}$, а также пользуясь весовым слабым типом $(1,1)$ у оператора $P_2^u$ (который, напомним, выводится из теоремы 4 статьи \cite{KisAnis}), имеем
	\begin{gather*}
		I_4
		\lesssim
		\int \sum 2^j \lambda_j(z_1)^{-\frac{1}{p}} \left( \int \abs*{g \psi_j^4} w(z_1, z_2) dz_2\right)^{\frac{1}{p}} v_j(z_1) dz_1
		+ \\ +
		\int \sum 2^j \left( \int \abs*{H(r_j (\gamma_j-1))^{k p}} dz_2 \right)^{\frac{1}{p}} dz_1
		\lesssim
		\dots
	\end{gather*}
	пользуясь ограниченностью оператора $H$ в решетке $L^k(\mathbb{T}, L^{kp}(\mathbb{T}))$, продолжим оценку
	\begin{gather*}
		\dots
		\lesssim
		\int \sum 2^j \lambda_j(z_1)^{-\frac{1}{p}} \left( \int \abs*{g \psi_j^4} w(z_1, z_2) dz_2\right)^{\frac{1}{p}} v_j(z_1) dz_1
		+ 
		\int \sum 2^j \left( \int \abs*{(\gamma_j-1) r_j}^{kp} dz_2 \right)^{\frac{1}{p}} dz_1
		\lesssim \\ \lesssim
		\int \int g(z_1, z_2) w(z_1, z_2) a(z_2) dz_2 dz_1
		+ 
		\int \sum 2^j \left( \int \abs*{(\gamma_j-1)}^{kp} w(z_1, z_2) dz_2 \right)^{\frac{1}{p}} v_j(z_1) dz_1
		\lesssim \\ \lesssim
		\int \int g(z_1, z_2) w(z_1, z_2) a(z_2) dz_2 dz_1
		+
		\int \sum 2^j v_j(z_1) \lambda_j(z_1)^{-\frac{1}{s}} \left( \int\limits_{\gamma_j > 1} \abs*{P_2^u g_1^j}^{\frac{p}{s}} w(z_1, z_2) dz_2 \right)^{\frac{1}{p}} dz_1.
	\end{gather*}
	Здесь первое слагаемое уже равняется $A$. Вскоре через $A$ оценим и второе. \\

	Для того, чтобы продолжить оценки, нам понадобится некоторое общее следствие из наличия у оператора слабого типа (1,1). Его можно найти, например, в \cite{stein} (8.15 в Гл. 1), но здесь мы приведем его в несколько измененной формулировке с доказательством. Пусть некоторый оператор $Q$ обладает слабым типом (1,1) относительно меры $\nu$, $A$ --- его (1,1)-норма; $\lambda > 0$, $0 < \alpha < 1$ --- некоторые числа; $e = \Set*{z}{\abs*{Qf}(z) > \lambda}$; $\rho_g(t)$ --- функция распределения функции $g$ по мере $\nu$ в точке $t$, тогда

	\begin{gather*}
	\int\limits_{e} \abs*{Qf}^\alpha(z) d\nu(z)
	=
	\int\limits_0^\infty t^{\alpha-1} \rho_{\chi_e \abs*{Q f}} (t) dt
	=
	\int\limits_0^c t^{\alpha-1} \rho_{\chi_e \abs*{Q f}} (t) dt + \int\limits_c^\infty t^{\alpha-1} \rho_{\chi_e \abs*{Q f}} (t) dt
	\leq \\ \leq
	\nu(e)\int\limits_0^c t^{\alpha-1}dt + A\norm*{f}\int\limits_c^\infty t^{\alpha-2}dt
	=
	\nu(e)\frac{1}{\alpha}c^{\alpha} + A \norm*{f} \frac{1}{1-\alpha}c^{\alpha-1}
	\end{gather*}
	Минимизируя последнее значение по $c$, получаем, что

	$$\int\limits_{e} \abs*{Qf}^\alpha(z) d\nu(z) \leq C(\alpha) A^\alpha \norm*{f}^\alpha \nu(e)^{1-\alpha}.$$

	Теперь, применяя только что доказанное общее утверждение, получим:
	\begin{gather*}
		\left( \int\limits_{\gamma_j > 1} \abs*{P_2^u g_1^j}^{\frac{p}{s}} w(z_1, z_2) dz_2 \right)^{\frac{1}{p}}
		\lesssim
		\left( \left((w(z_1, \cdot) \mu) \left( \Set*{z_2}{\gamma_j(z_1, z_2) > 1}\right)\right)^{\frac{s-p}{p}} \int \abs*{g \psi_j^4} w(z_1, z_2) dz_2\right)^{\frac{1}{s}}
		\lesssim \\ \lesssim
		\left( \lambda_j(z_1)^{-1} \int \abs*{g \psi_j^4} w(z_1, z_2) dz_2 \right)^{\frac{1}{p} - \frac{1}{s}} \left( \int \abs*{g \psi_j^4} w(z_1, z_2) dz_2 \right)^{\frac{1}{s}}
		=
		\lambda_j(z_1)^{\frac{1}{s}-\frac{1}{p}} \left( \int \abs*{g \psi_j^4} w(z_1, z_2) dz_2 \right)^{\frac{1}{p}}
	\end{gather*}

	Заметим, что константы, скрытые выше под символом ``$\lesssim$'' не зависят от переменной $z_1$, так как в силу равномерности наложенных условий Макенхаупта, $(1,1)$-норма проектора $P_2^u$, зависящая (в силу теоремы 4 из \cite{KisAnis}) только от соответствующих $A_p$-норм весов $w$ и $w_1$, может быть оценена константой, не зависящей от $z_1$. \\

	Так что

	$$I_4 \lesssim A + \int \sum 2^j v_j(z_1) \lambda_j(z_1)^{-\frac{1}{p}} \left( \int \abs*{g \psi_j^4} w(z_1, z_2) dz_2 \right)^{\frac{1}{p}} dz_1 \lesssim A,$$

	что и требовалось доказать.

\end{proof}
						\subsection*{\begin{center}Возвращаемся к пространствам Харди\end{center}}

Теперь, наконец, можно получить утверждение про $K$-замкнутость пары
$$\left(H_p(b_1(z_1) w_1(z_1, z_2) a_1(z_2)), H_\infty(b_2(z_1) w_2(z_1,z_2) a_2(z_2))\right)$$
в паре
$$\left(L_p(b_1(z_1) w_1(z_1, z_2) a_1(z_2)), L_\infty(b_2(z_1) w_2(z_1,z_2) a_2(z_2))\right).$$
Верна следующая теорема:
\begin{theorem}
	\label{inf_neib_all_q}
	K-замкнутость имеет место, если
	\begin{enumerate}[label={\arabic*)}]
		\item[0)] $w_1, w_2, w_1^{\frac{1}{1-p}} \in L^1(\mathbb{T}^2)$
		\item $w_2^p w_1$ удовлетворяет условию $A_{\infty}$ по второй переменной равномерно,
		\item $w_1$ удовлетворяет условию $A_p$ по второй переменной равномерно,
		\item $w_2$ удовлетворяет условию $A_1$ по второй переменной равномерно,
		\item $\log(a_i), \log(b_i) \in \text{BMO}$,
		\item $\log(w_1(\cdot, z_2))$, $\log(w_2(\cdot, z_2))$ лежат в пространстве $\text{BMO}$ по первой переменной равномерно.
	\end{enumerate}
\end{theorem}
\begin{remark}
	Заменив одномерные $A_p$-условия на соответствующие им двумерные, мы, хоть и несколько ослабив предыдущие две теоремы, избавимся от многих неприятных технических условий и сделаем теоремы гораздо более удобными для использования. Например, пропадут условия \textit{0)} в обеих теоремах (так как они будут следовать из двумерных $A_p$). Окончательный результат, получающийся при таком ослаблении теорем, был изложен в теореме \ref{rght} во введении.
\end{remark}

\newpage
					\section*{\begin{center}\S4 $K$-замкнутость при $r = 1$, $p = \infty$ (склейка)\end{center}}
						В этой части мы хотим получить аналог предложения 7 из \cite{ruts}, а именно теорему \ref{glue} из введения. Напомним ее формулировку. 
\begingroup
\def\thetheorem{\ref{glue}}
\begin{theorem}
	Если $w_1, w_2 \in A_1$ и $w_1 w_2 \in A_\infty$ (оба условия Макенхаупта двумерные), то пара 
	$$\left(H_1(w_1(z_1, z_2)), H_\infty(w_2(z_1,z_2))\right)$$
	$K$-замкнута в паре
	$$\left(L_1(w_1(z_1, z_2)), L_\infty(w_2(z_1,z_2))\right).$$
\end{theorem}
\addtocounter{theorem}{-1}
\endgroup
\begin{proof}
	Нам достаточно найти такие числа $0 < \theta_1 < \theta_2 < \theta_3 < \theta_4 < 1$, что $K$-замкнутость соответствующих пар пространств Харди имеет место для пар (про интерполяцию взвешенных пространств Лебега см. \cite{Fre})
	$$\left(L_1(w_1), L_{\frac{1}{1-\theta_2}}(w_1 w_2^{\frac{\theta_2}{\theta_2-1}})\right),$$
	$$\left(L_{\frac{1}{1-\theta_1}}(w_1 w_2^{\frac{\theta_1}{\theta_1-1}}), L_{\frac{1}{1-\theta_4}}(w_1 w_2^{\frac{\theta_4}{\theta_4-1}})\right),$$
	$$\left(L_{\frac{1}{1-\theta_3}}(w_1 w_2^{\frac{\theta_3}{\theta_3-1}}), L_\infty(w_2)\right),$$
	остальное сделает теорема типа Вольфа для $K$-замкнутости из \cite{KisXu}. \\
	
	Выберем произвольные $0 < \theta_1 < \theta_2 < \theta_3 < \theta_4 < 1$. Действительно, из факторизационной теоремы Джонса (см. в \cite{stein}) следует, что $w_1 w_2^{\frac{\theta}{\theta-1}} \in A_{\frac{1}{1-\theta}}$, что сразу же дает $K$-замкнутость для пары $\left(L_{\frac{1}{1-\theta_1}}(w_1 w_2^{\frac{\theta_1}{\theta_1-1}}), L_{\frac{1}{1-\theta_4}}(w_1 w_2^{\frac{\theta_4}{\theta_4-1}})\right)$, а в сумме с теоремой \ref{one_naib} дает $K$-замкнутость еще и для пары $\left(L_1(w_1), L_{\frac{1}{1-\theta_2}}(w_1 w_2^{\frac{\theta_2}{\theta_2-1}})\right)$. Выполнено, по условию, что $w_2^{\frac{1}{1-\theta_3}} w_1 w_2^{\frac{\theta_3}{\theta_3-1}} = w_1 w_2 \in A_{\infty}$, кроме того, из факторизационной теоремы Джонса снова выводим, что $w_1 w_2^{\frac{\theta_3}{\theta_3-1}}$ удовлетворяет условию $A_{\frac{1}{1-\theta_3}}$, а значит теорема \ref{inf_neib_all_q} дает $K$-замкнутость для пары $\left(L_{\frac{1}{1-\theta_3}}(w_1 w_2^{\frac{\theta_3}{\theta_3-1}}), L_\infty(w_2)\right)$. Доказательство закончено.

\end{proof}

					\section*{\begin{center} Список литературы \end{center}}


\begin{thebibliography}{9}

	\bibitem{ruts}
		Д. В. Руцкий,
		\emph{Весовое разложение Кальдерона–Зигмунда и некоторые его приложения к интерполяции},
		Зап. научн. сем. ПОМИ, 2014, том 424, 186–200,
		2014

	\bibitem{xu}
		Quanhua Xu,
		\emph{Some properties of the quotient space $(L^{1}(\mathbf{T}^{d})/H^{1}(D^{d}))$},
		Illinois J. Math. Volume 37, Issue 3 (1993), 437-454,
		1993

	\bibitem{uncond}
		J. García-Cuerva, K. Kazarian,
		\emph{Calderón-Zygmund operators and unconditional bases of weighted Hardy spaces},
		Studia Mathematica Volume 109, Issue 3, 255-276,
		1994
	
	\bibitem{Isra}
		D. Israfilov, A. Guven,
		\emph{Approximation by trigonometric polynomials in weighted Orlicz spaces},
		Studia Mathematica 174 (2006), 147-168,
		2006

	\bibitem{kis}
		S. V. Kislyakov,
		\emph{Interpolation Involving Bounded Bianalytic Functions},
		Operator Theory: Advances and Applications, Vol. 113,
		2000

	\bibitem{kis2}
		S. V. Kislyakov,
		\emph{Interpolation of $H^p$ spaces: some recent developments}
		Israel Math. Conf. Proceedings, 13 (1999), 102-140,
		1999

	\bibitem{KisAnis}
		S. V. Kislyakov, D. S. Anisimov,
		\emph{Double singular integrals: interpolation and correction},
		St. Petersburg Mathematical Journal, 2005, 16:5, 749–772,
		2005

	\bibitem{GamKis}
		T. W. Gamelin, S. V. Kislyakov,
		\emph{Chapter 16 - Uniform Algebras as Banach Spaces},
		Handbook of the Geometry of Banach Spaces, Vol. 1, 671-706,
		2001

	\bibitem{KisBurgProj}
		S. V. Kislyakov,
		\emph{Bourgain's Analytic Projection Revisited},
		Proceedings of the American Mathematical Society, Vol. 126, No. 11 (Nov., 1998), pp. 3307-3314,
		1998

	\bibitem{KisAbsSumm}
		С. В. Кисляков,
		\emph{Абсолютно суммирующие операторы на диск-алгебре},
		Алгебра и анализ, 1991,	том 3, выпуск 4, страницы 1–77,
		1991

	\bibitem{KisXu}
		С. В. Кисляков, Куанхуа Шу,
		\emph{Вещественная интерполяция и сингулярные интегралы},
		Алгебра и анализ, 1996,	том 8, выпуск 4, страницы 75–109,
		1996

	\bibitem{Fre}
		D. Freitag,
		\emph{Real Interpolation of Weighted $L_p$-Spaces},
		Mathematische Nachrichten, Vol. 86, Issue 1 (1978), 15-18,
		1978

	\bibitem{coifmanweiss}
		R. R. Coifman, G. Weiss,
		\emph{Extensions of Hardy spaces and their use in analysis},
		Bull. Amer. Math. Soc. Volume 83, Number 4 (1977), 569-645,
		1977

	\bibitem{stein}
		E. M. Stein,
		\emph{Harmonic Analysis: Real-Variable Methods, Orthogonality, and Oscillatory Integrals},
		Princeton University Press,
		1993
	
	\bibitem{Rud}
		W. Rudin,
		\emph{Function theory in polydiscs},
		Mathematics lecture note series (Выпуск 41), W. A. Benjamin,
		1969

	\bibitem{Rubio}
		J. García-Cuerva, J.L. Rubio de Francia,
		\emph{Weighted Norm Inequalities and Related Topics},
		North-Holland mathematics studies 116,
		1985

	\bibitem{Strom}
		J. O. Stromberg, A. Torchinsky,
		\emph{Weighted Hardy Spaces},
		Lecture Notes in Mathematics, Springer-Verlag,
		1989

	\bibitem{Woj}
		P. Wojtaszczyk,
		\emph{Banach Spaces for Analysts},
		Cambridge University Press,
		1991

	\bibitem{Bergh}
		J. Bergh, J. Lofstrom,
		\emph{Interpolation Spaces},
		Springer Berlin Heidelberg,
		1976
		
	\bibitem{KantAk}
		Л. В. Канторович, Г. П. Акилов,
		\emph{Функциональный анализ},
		Москва ``Наука'',
		1984

\end{thebibliography}
\end{document}